\documentclass[11pt]{article}
\usepackage{a4wide}

\usepackage[style=numeric,sorting=none]{biblatex}
\usepackage[american]{babel}
\usepackage[fleqn]{amsmath}
\usepackage{amssymb}
\usepackage{tabularx}
\usepackage{ltablex}
\usepackage{csquotes}
\usepackage{graphicx}
\usepackage{placeins}
\usepackage{eurosym}
\keepXColumns
\usepackage[dvipsnames]{xcolor}
\usepackage{array} 

\newcolumntype{C}[1]{>{\centering\arraybackslash}m{#1}}

\title{Multi-year Investment Modelling in Energy Systems}
\author{Diego A. Tejada-Arango}
\date{July 2023}

\begin{document}

\maketitle

\section{Nomenclature}

\subsection{Indices}
\begin{tabularx}{\textwidth}{| l | X |} 
\hline 
\textbf{Name} & \textbf{Description} \\ 
\hline 
\endhead 
$y,m,\mu$ & years                  \\ 
$k$       & representative periods \\ 
$t$       & time steps             \\ 
\hline 
\end{tabularx}

\subsection{Sets}
\begin{tabularx}{\textwidth}{| l | X |} 
\hline 
\textbf{Name} & \textbf{Description} \\ 
\hline 
\endhead 
$\mathcal{M}  $  & Milestone years      \\ 
\hline 
\end{tabularx}

\subsection{Parameters}
\begin{tabularx}{\textwidth}{| l | X | l |} 
\hline 
\textbf{Name} & \textbf{Description} & \textbf{Units} \\ 
\hline 
\endhead 
$R$           & Social Discount Rate      & [\%] \\ 
$WACC_{y}$    & Weighted Average Cost of Capital at year $y$& [\%] \\
$LT$          & Lifetime of the asset  & [years] \\
$Y$           & Total number of years in the time horizon & [years] \\ 
$C_{y}^{T}$   & Total investment cost at year $y$ (i.e., overnight cost) & [\euro/MW] \\
$C_{y}^{A}$   & Annualized investment cost  at year $y$ & [\euro/MW/year] \\
$C^{I}$       & Total investment cost & [\euro] \\
$C^{O}$       & Total operational cost & [\euro] \\
$SV_{y}$      & Salvage value & [\euro] \\
$SVP_{y}$     & Salvage value percentage & [\%] \\
$C_{y}^{op}$  & Operational cost & [\euro/MWh] \\
$W_{yk}^{op}$ & Operational cost weight of representative $k$ in year $y$ & [h] \\
$W_{y}^{I}$   & Milestone year weight & [years] \\
$W_{my}^{O}$  & Weight of milestone year $m$ on year $y$ & [0-1] \\
$W_{m\mu}^{M}$  & Weight of milestone year $m$ on operational year $\mu$ & [0-1] \\
$W_{my\mu}^{Y}$ & Weight of milestone year $m$ on year $y$ and operational year $\mu$ & [0-1] \\
\hline 
\end{tabularx}

\subsection{Variables}
\begin{tabularx}{\textwidth}{| l | X | l |} 
\hline 
\textbf{Name} & \textbf{Description} & \textbf{Units}  \\ 
\hline 
\endhead 
$x_{y}$    & capacity investment at year $y$   & [MW]   \\ 
$p_{ykt}$  & production at period $t$, representative $k$, and year $y$   & [MW]   \\ 
$p_{mykt}$ & production at period $t$, representative $k$, operational year $y$, from the investment made in year $m$   & [MW]   \\
\hline 
\end{tabularx}

\section{Introduction} \label{sec:intro}
This document summarises the main multi-year investment modelling approaches in energy planning models. Therefore, here we will go from a simple (basic) formulation to a more complex (general) one to understand different levels of detail, including examples to make more accessible the understanding of the concepts.

In the literature, there are two approaches to managing multi-horizon investments:
\begin{itemize}
    \item \textit{Total Cost}: Using the total cost (a.k.a. overnight cost) for the investment year.

\begin{equation}
 \label{eq:total_cost_from_annualized_cost}
    C_{y}^{T} = \sum_{j = y}^{y+LT-1}\frac{1}{(1+WACC_{y})^{j - y}} \cdot C_{j}^{A} \quad \forall y
\end{equation}

    \item \textit{Annualized Cost}: Distributing the total cost over the years the asset is active, which includes all the years following the build year until its lifetime.

\begin{equation}
 \label{eq:annualized_cost_from_total_cost}
    C_{y}^{A} = \frac{WACC_{y}}{(1+WACC_{y}) \cdot \left( 1- \frac{1}{(1+WACC_{y})^{LT}} \right)} \cdot C_{y}^{T} \quad \forall y
\end{equation}

\end{itemize}

It's important to note that equation (\ref{eq:total_cost_from_annualized_cost}) has infinite solutions for defining the annuities ($C_{j}^{A}$). One standard method is to assume that the annuities remain constant throughout the lifetime, as illustrated in equation (\ref{eq:annualized_cost_from_total_cost}). In addition, equation (\ref{eq:annualized_cost_from_total_cost}) does not discount the first year of the investment, which is explained in \cite{Centeno2018}. Another common approach is discounting the first year and using equation (\ref{eq:annualized_cost_from_total_cost_alternative}) instead. Depending on the assumptions made, either equation (\ref{eq:annualized_cost_from_total_cost}) or (\ref{eq:annualized_cost_from_total_cost_alternative}) can be used to calculate the annualized cost from the total cost.

\begin{equation}
 \label{eq:annualized_cost_from_total_cost_alternative}
    C_{y}^{A} = \frac{WACC_{y}}{1- \frac{1}{(1+WACC_{y})^{LT}}} \cdot C_{y}^{T} \quad \forall y
\end{equation}

Let us consider the following data as an example for the first year ($y=0$):

\begin{itemize}
    \item Total capital cost $C_{0}^{T}=100$ [\euro/MW]
    \item The $WACC_{0}=2$ \%
    \item The economic lifetime $LT=5$ years
\end{itemize}

Therefore, the annualized cost for each case is:

\begin{itemize}
    \item Annualized capital cost using equation (\ref{eq:annualized_cost_from_total_cost}): $C_{0}^{I}=20.80$ [\euro/MW]
    \item Annualized capital cost using equation (\ref{eq:annualized_cost_from_total_cost_alternative}): $C_{0}^{I}=25.05$ [\euro/MW]
\end{itemize}

Although the values seem different, they will represent the total investment cost if correctly discounted for each case. For instance, using the annualized value from equation (\ref{eq:annualized_cost_from_total_cost}):

\begin{equation*}
    C_{0}^{T} = 100 = 20.80 \cdot \left(1+\frac{1}{(1+0.02)^{1}}+\frac{1}{(1+0.02)^{2}}+\frac{1}{(1+0.02)^{3}}+\frac{1}{(1+0.02)^{4}} \right)
\end{equation*}

Whereas using the annualized value from equation (\ref{eq:annualized_cost_from_total_cost_alternative}) the calculation is:

\begin{equation*}
    C_{0}^{T} = 100 = 25.05 \cdot \left(\frac{1}{(1+0.02)^{1}}+\frac{1}{(1+0.02)^{2}}+\frac{1}{(1+0.02)^{3}}+\frac{1}{(1+0.02)^{4}} +\frac{1}{(1+0.02)^{5}}\right)
\end{equation*}

Therefore, using equation (\ref{eq:annualized_cost_from_total_cost_alternative}) means that equation (\ref{eq:total_cost_from_annualized_cost}) changes slightly to represent the total cost adequately. For the remainder of this document, we will utilize equations (\ref{eq:total_cost_from_annualized_cost}) and (\ref{eq:annualized_cost_from_total_cost}) to continue our discussion.

The reference \cite{Brown2020} provides a concise summary of the approach utilized in popular planning models within the energy sector, including TIMES, DIMENSION, SWITCH, and PRIMES. In addition, it suggests that the methods are only comparable under the condition that $Y \rightarrow \infty$. Otherwise, this could lead to different interpretations of the costs in the objective function. However, this report demonstrates that with careful execution, the objective function can be made identical in both techniques.

\section{Basic formulations}
Here we will go from a simple (basic) formulation to a more complex (detailed) one for each method. Equation (\ref{eq:objective_function}) shows a general formulation of the objective function for a multi-year investment planning model, while equation (\ref{eq:production_limit}) shows the relationship between the variables where the investment decisions limit generation production considering its lifetime.

\begin{equation}
 \label{eq:objective_function}
    \min_{x_{y},p_{ykt}} \quad C^{I} + C^{O} 
\end{equation}

\begin{equation}
 \label{eq:production_limit}
 \begin{aligned}
    &s.t. \\
    &p_{ykt} \leq \sum_{j=\max(y-LT+1,0)}^{y}x_{j} \quad \forall y \forall k \forall t
 \end{aligned}
\end{equation}

Where $C^{I}$ refers to the investment cost and $C^{O}$ to the operational cost. A general definition for $C^{O}$ is shown in equation (\ref{eq:basic_operational_cost_formulation}):

\begin{equation}
 \label{eq:basic_operational_cost_formulation}
    C^{O} = \sum_{y=0}^{Y-1}\frac{1}{(1+R)^{y}} \cdot C_{y}^{op} \sum_{k}W_{yk}^{op} \cdot \sum_{t} p_{ykt} \\
\end{equation}

As for the investment cost component in the objective function $C^{I}$, this document will focus on obtaining this component in both approaches: the total and annualized investment costs.
There are additional costs that have not been included in this objective function for the sake of simplicity, such as fixed costs, decommissioning costs, and costs for mothballing.

\subsection{Total Investment Cost}
Equation (\ref{eq:basic_total_cost_formulation}) shows the total investment cost through the modelled years, considering the social discount rate $R$.

\begin{equation}
 \label{eq:basic_total_cost_formulation}
 \begin{aligned}
    & C^{I} = \sum_{y=0}^{Y-1}\frac{1}{(1+R)^{y}} \cdot C_{y}^{T}x_{y} \\
    & C^{I} = C_{0}^{T}x_{0} + \frac{1}{(1+R)^{1}} \cdot C_{1}^{T}x_{1} + \frac{1}{(1+R)^{2}} \cdot C_{2}^{T}x_{2} + \cdots
 \end{aligned}
\end{equation}

\subsection{Annualized Investment Cost}
Equation (\ref{eq:basic_annualized_cost_formulation}) is obtained by replacing (\ref{eq:total_cost_from_annualized_cost}) in (\ref{eq:basic_total_cost_formulation}). Bear in mind that the value of the parameter $C_{y}^{A}$ is determined by (\ref{eq:annualized_cost_from_total_cost}).

\begin{equation}
 \label{eq:basic_annualized_cost_formulation}
 \begin{aligned}
    C^{I} = & \sum_{y=0}^{Y-1}\frac{1}{(1+R)^{y}} \cdot x_{y} \cdot \sum_{j = y}^{\min(y+LT-1,Y-1)}\frac{1}{(1+WACC_{y})^{j - y}} \cdot C_{y}^{A} \\
    C^{I} = & C_{0}^{A} \cdot \left(1 + \frac{1}{(1+WACC_{0})^{1}} + \frac{1}{(1+WACC_{0})^{2}} + \cdots \right) \cdot x_{0}\\ 
            & + \frac{1}{(1+R)^{1}} \cdot C_{1}^{A} \cdot \left( 1 + \frac{1}{(1+WACC_{1})^{1}} + \frac{1}{(1+WACC_{1})^{2}} + \cdots \right) \cdot x_{1} \\ 
            & + \frac{1}{(1+R)^{2}} \cdot C_{2}^{A} \cdot \left( 1 + \frac{1}{(1+WACC_{2})^{1}} + \frac{1}{(1+WACC_{2})^{2}} + \cdots \right) \cdot x_{2} + \cdots
 \end{aligned}
\end{equation}

\section{Analysis}
As mentioned in Section \ref{sec:intro}, equations (\ref{eq:basic_total_cost_formulation}) and (\ref{eq:basic_annualized_cost_formulation}) are equivalent when considering yearly investment and $Y \rightarrow \infty$. This statement can be checked by considering the following relationships:

\begin{equation*}
 \label{eq:equivalence_total_and_annualized_cost}
 \begin{aligned}
    C_{0}^{T} = & C_{0}^{A} \cdot \left(1 + \frac{1}{(1+WACC_{0})^{1}} + \frac{1}{(1+WACC_{0})^{2}} + \cdots \right) \\ 
    C_{1}^{T} = & C_{1}^{A} \cdot \left(1 + \frac{1}{(1+WACC_{1})^{1}} + \frac{1}{(1+WACC_{1})^{2}} + \cdots \right) \\ 
    C_{2}^{T} = & C_{2}^{A} \cdot \left(1 + \frac{1}{(1+WACC_{2})^{1}} + \frac{1}{(1+WACC_{2})^{2}} + \cdots \right) \\
    \vdots \qquad & \vdots
 \end{aligned}
\end{equation*}

It is also essential to highlight that no matter the values of the social discount rate ($R$) or the technology-specific discount rate ($WACC$), the total investment cost ($C^{I}$) in the objective function is the same. This affirmation holds if you obtain the annualized value using equation (\ref{eq:annualized_cost_from_total_cost}) and the same $WACC$ value in the equation (\ref{eq:basic_annualized_cost_formulation}). If this isn't the situation, there's a possibility of having different total values in the objective function. In that situation, unifying the $WACC$ and using the same value for both equations is better to avoid errors and confusion.

Finally, the $Y \rightarrow \infty$ assumption is not realistic to apply. The usual use case considers that modelled years do not cover the entire investment lifetime and use milestone years (or representative years) to reduce the computational burden. These two cases are discussed in the following sections.

\subsection{Having Modelled Years $<$ Lifetime}
First, Let's analyze the case of years that do not cover the entire investment lifetime but still use yearly investment decisions. For instance, if the last modelled year ($Y$) is 2050 and we have an investment made in 2040 with a lifetime of 15 years, then the last four years will not be considered in the optimisation. In this situation:

\begin{itemize}
    \item The annualized cost approach in (\ref{eq:basic_annualized_cost_formulation}) will only recover the annualized investment costs from 2040 to 2050, and therefore it assumes that the non-modelled years (i.e., 2051 to 2054) will be recovered independently from the modelled years and are not part of the optimization. Under this assumption, there is no need for changing the equation (\ref{eq:basic_annualized_cost_formulation}) due to modelling years below the asset's lifetime.
    \item The total investment cost approach has a problem using the equation (\ref{eq:basic_total_cost_formulation}) in the objective function. It will consider the total cost of the investment but will miss the operational costs from the non-modelled years to compensate for the total investment cost. In our example, the operational cost from 2040 to 2050 must compensate for the total investment cost of a 15-year lifetime asset. To avoid this situation, the standard approach is to consider the salvage value ($SV$).
\end{itemize}

\subsubsection{The Salvage value}
The salvage value is what a company expects to get by selling or disassembling an asset at the end of its useful life, which will be the last modelled year in our case for these types of models. As in the annualized cost approach, one common assumption is that the salvage value is equivalent to the cost outside the modelled years (i.e., the non-modelled years will be recovered independently of the optimization problem). Therefore, the salvage value that satisfies this assumption can be determined from the annualized cost value, see equation (\ref{eq:salvage_value}).

\begin{equation}
 \label{eq:salvage_value}
    SV_{y} = C_{y}^{A} \cdot \sum_{j = Y+1}^{y+LT-1}\frac{1}{(1+WACC_{y})^{j - y}} \quad \forall y \in [0,Y]
\end{equation}

The salvage value $SV$ will depend on the year $y$, the last modelled year $Y$, the annualized value $C_{y}^{A}$, the lifetime $LT$, and the $WACC$. For instance, consider the following data for the first year ($y=0$):

\begin{itemize}
    \item Total capital cost $C_{0}^{T}=100$ [\euro/MW]
    \item The $WACC_{0}=5$ \%
    \item The last modelled year $Y=4$
    \item The economic lifetime $LT=8$ years
    \item Annualized capital cost using equation (\ref{eq:annualized_cost_from_total_cost}): $C_{0}^{I}=14.74$ [\euro/MW]
\end{itemize}

Therefore, using equation (\ref{eq:salvage_value}), the $SV_{0} = 33.01$ [\euro/MW] corresponds to the remaining cost in non-modelled years: 5, 6, and 7. The $SV$ is referenced to the first year $y=0$ and can be considered as a benefit of the investment in the objective function, as in equations (\ref{eq:total_cost_formulation_with_salvage_value_a}) and (\ref{eq:total_cost_formulation_with_salvage_value_b}).

\begin{equation}
 \label{eq:total_cost_formulation_with_salvage_value_a}
    C^{I} = \sum_{y=0}^{Y-1}\frac{1}{(1+R)^{y}} \cdot (C_{y}^{T}-SV_{y})\cdot x_{y}
\end{equation}

\begin{equation}
 \label{eq:total_cost_formulation_with_salvage_value_b}
    C^{I} = \sum_{y=0}^{Y-1}\frac{1}{(1+R)^{y}} \cdot C_{y}^{T}(1-SVP_{y})\cdot x_{y} 
\end{equation}

Where, 
\begin{equation*}
    SVP_{y}=\frac{SV_{y}}{C_{y}^{T}}    \quad   \forall y
\end{equation*}

Equations (\ref{eq:total_cost_formulation_with_salvage_value_a}) and (\ref{eq:total_cost_formulation_with_salvage_value_b}) are equivalent and can be used instead of equation (\ref{eq:basic_total_cost_formulation}) to account for the salvage value in the total investment cost for the objective function.

If the $SV$ is obtained using the same assumptions as the ones for the annualized cost method, then considering the $SV$ will make both methods equivalent. In the final year of modelling, $Y$, the $SV$ will make up for LT-1 years of the asset. This implies that the total investment method will only consider the cost of one year in the last modelled year, resulting in the same value as the annualized cost approach. As a result, both methods represent the same cost in the objective function.

\subsection{Using Milestone Years} \label{sec:milestone_years_weight}
When studying pathways to the year 2050, it is impractical to model yearly investment decisions on a large scale. Therefore, a common approach is to use the so-called milestone years, where the major investment decisions will be made (e.g., years 2030, 2035, 2040, and 2050). Each milestone year represents part of the non-modelled years.

Consider the information in Table \ref{tab:milestone_years_ex1} as an example; therefore, years 1, 3, and 5 are not modelled since the milestone years for investments are years 0, 2, and 5. 

\begin{table}[ht!]
\centering
\begin{tabular}{|c|c|c|c|c|c|c|}
\hline
\textbf{Years $y$} & 0 & 1 & 2 & 3 & 4 & 5\\ 
\hline
\hline
\textbf{Milestone year $m$} & \textcolor{blue}{0} & - & \textcolor{blue}{2} & - & - & \textcolor{blue}{5}\\ 
\hline
\textbf{Milestone year weight $W_{m}^{I}$} & \textcolor{blue}{2} & - & \textcolor{blue}{3} & - & - & \textcolor{blue}{1} \\ 
\hline
\end{tabular}
\caption{Milestone years}
\label{tab:milestone_years_ex1}
\end{table}

The \textit{total investment cost} approach needs a weight coefficient in the operational cost $W_{y}^{I}$ to properly account for the non-modelled years. In this example, year 0 has a weight of 2, year 2 has a weight of 3, whereas year 4 equals 1. Moreover, since the total investment cost is considered in the milestone year, the total investment cost $C^{I}$ is not needed to multiply by the coefficient $W_{y}^{I}$. Equation (\ref{eq:total_investment_cost_milestone_weights}) shows the formulation, including all these considerations:

\begin{equation}
 \label{eq:total_investment_cost_milestone_weights}
  \begin{aligned}
    & \min_{x_{m},p_{mkt}} \quad C^{I} + C^{O}  \\
    & C^{I} = \sum_{m \in \mathcal{M}}\frac{1}{(1+R)^{m}} \cdot (C_{m}^{T}-SV_{m})\cdot x_{m} \\ 
    & C^{O} = \sum_{m \in \mathcal{M}}\frac{1}{(1+R)^{m}} \cdot \textcolor{blue}{W_{m}^{I}} \cdot C_{m}^{op} \sum_{k}W_{mk}^{op} \cdot \sum_{t} p_{mkt}
  \end{aligned}
\end{equation}    

The \textit{annualized investment cost} approach does not need a weight coefficient a priory in its formulation if one assumes that the not modelled years recover their own cost (i.e., the same assumption as for the years after the last modelled year). However, it might be desirable to include the weight coefficient $W_{m}^{I}$ for this method if the milestone years in unevenly spaced and one wants to consider that a particular milestone year has more weight than the others. For instance, the asset has more production costs at the beginning and less at the end of the horizon, such as the production costs at the beginning compensate enough for investment cost. Therefore, using the weight coefficient $W_{m}^{I}$ in this method is also usual practice; see equation (\ref{eq:annualized_investment_cost_milestone_weights}).

\begin{equation}
 \label{eq:annualized_investment_cost_milestone_weights}
  \begin{aligned}
    & \min_{x_{m},p_{mkt}} \quad C^{I} + C^{O}  \\
    & C^{I} = \sum_{m \in \mathcal{M}}\frac{1}{(1+R)^{m}} \cdot x_{m} \cdot \sum_{j = m | j \in \mathcal{M}}^{\min(m+LT-1,Y-1)}\frac{1}{(1+WACC_{m})^{j - m}} \cdot \textcolor{blue}{W_{j}^{I}} \cdot C_{m}^{A} \\ 
    & C^{O} = \sum_{m \in \mathcal{M}}\frac{1}{(1+R)^{m}} \cdot \textcolor{blue}{W_{m}^{I}} \cdot C_{m}^{op} \sum_{k}W_{mk}^{op} \cdot \sum_{t} p_{mkt}
  \end{aligned}
\end{equation} 

Please take note that the calculation for investment cost in equation (\ref{eq:annualized_investment_cost_milestone_weights}) only considers milestone years in the inner summation, owing to $j \in \mathcal{M}$.

Let us continue with the example in Table \ref{tab:milestone_years_ex1} (i.e., $\mathcal{M} = \{0,2,5\}$ and  $W_{m}^{I} = [2,3,1]$) and assume that $LT=6$. Therefore, the investment cost in equations in the previous equations becomes:

\textit{Equation (\ref{eq:total_investment_cost_milestone_weights}) - total investment cost method:}
\begin{equation*}
    C^{I} = \textcolor{orange}{C_{0}^{T}}x_{0} + \frac{1}{(1+R)^{2}} (C_{2}^{T}-SV_{2})\cdot x_{2} + \frac{1}{(1+R)^{5}} (C_{5}^{T}-SV_{5})\cdot x_{5}
\end{equation*}    

\textit{Equation (\ref{eq:annualized_investment_cost_milestone_weights}) - annualized investment cost method:}
\begin{equation*}
  \begin{aligned}
     C^{I} = & \textcolor{orange}{C_{0}^{A}\left(2+\frac{3}{(1+WACC_{0})^{2}} +\frac{1}{(1+WACC_{0})^{5}}\right)}\cdot x_{0} \\ 
             & +\frac{1}{(1+R)^{2}} \cdot C_{2}^{A} \left(3+\frac{1}{(1+WACC_{2})^{3}} \right) \cdot x_{2} + \frac{1}{(1+R)^{5}}C_{5}^{A\cdot}x_{5} \\             
  \end{aligned}
\end{equation*}    

By simply inspecting the methods used, it is clear that the costs represented in the objective function are not the same as before when yearly investment decisions were made. For example, if we compare the first term in both costs (highlighted in \textcolor{orange}{orange}), we can see that the values differ according to equations (\ref{eq:total_cost_from_annualized_cost}) and (\ref{eq:annualized_cost_from_total_cost}).

The operational cost is the same in both methods since the assumption here is that the milestone years represent the non-modelled years using the weight $W_{m}^{I}$. Applying the values of the example in equations (\ref{eq:total_investment_cost_milestone_weights}) and (\ref{eq:annualized_investment_cost_milestone_weights}), the operational cost is:

\begin{equation*}
  \begin{aligned}
    C^{O} = & \sum_{m \in \mathcal{M}}\frac{1}{(1+R)^{m}} \cdot \textcolor{blue}{W_{m}^{I}} \cdot C_{m}^{op} \sum_{k}W_{mk}^{op} \cdot \sum_{t} p_{mkt} \\
    & + 2 \cdot C_{0}^{op} \sum_{k}W_{0,k}^{op} \cdot \sum_{t} p_{0,kt} \\
    & + \frac{3}{(1+R)^{2}} \cdot C_{2}^{op} \sum_{k}W_{2,k}^{op} \cdot \sum_{t} p_{2,kt} \\
    & + \frac{1}{(1+R)^{5}} \cdot C_{5}^{op} \sum_{k}W_{5,k}^{op} \cdot \sum_{t} p_{5,kt} \\
  \end{aligned}
\end{equation*}

Therefore, this example assumes that milestone year 2 represents non-modelled years 3 and 4. The issue here is that years 3 and 4 should have a lower value of money over time than year 2, meaning that this approach overestimates the operational cost. Furthermore, this approach presents another issue. If the asset's lifetime ends between two milestone years, the $C^I$ and $C^O$ in the objective function cannot account for this and will assume the asset is available throughout the entire range between the milestone years. For instance, if $LT=4$ in Table (\ref{tab:milestone_years_ex1}), the asset should be available until year 3. However, equations (\ref{eq:annualized_cost_from_total_cost}) will include them in the objective function until the final year (i.e., $y=5$), resulting in an additional two years in the objective function (i.e., overestimation of the operational costs again).

\section{Proposals to Improve the Milestone Years Method} \label{sec:proposals_main_sec}
As discussed in the previous section, multiplying the milestone years by the weight $W_{m}^{I}$ will lead to two main issues in the objective function of both methods:

\begin{itemize}
    \item The value of money over time is not correctly represented for the non-modelled years.
    \item The lifetime of the assets is overestimated if the last year is between two milestone years.
\end{itemize}

In the following sections, we propose some formulations that attempt to solve these shortcomings in both approaches.

\subsection{Proposal 1 - Mapping investment milestone years to all years ($W_{my}^{O}$)} \label{sec:proposal_mapping_m_to_y}
The first proposal is based on adding an extra index $y$ in the objective function that includes all the years in the time  horizon (i.e., milestone and non-modelled years). Then, a mapping between the non-modelled years and the milestone year is needed to properly account for them. The weight of milestone year $m$ on year $y$ ($W_{my}^{O}$) ensures this mapping among years. For the total investment cost, we only need to change the operational cost in the objective function as shown in equation (\ref{eq:total_investment_cost_milestone_weights_proposal1}). The very change applies to the annualized investment cost approach, which will account additionally for all the years during the asset's lifetime in the $C^{I}$; see equation (\ref{eq:annualized_investment_cost_milestone_weights_proposal1}).

\textit{Total investment cost method:}
\begin{equation}
 \label{eq:total_investment_cost_milestone_weights_proposal1}
  \begin{aligned}
    & \min_{x_{m},p_{mkt}} \quad C^{I} + C^{O}  \\
    & C^{I} = \sum_{m \in \mathcal{M}}\frac{1}{(1+R)^{m}} \cdot (C_{m}^{T}-SV_{m})\cdot x_{m} \\ 
    & C^{O} = \sum_{y=0}^{Y-1}\frac{1}{(1+R)^{y}}\sum_{m \in \mathcal{M}} \textcolor{blue}{W_{my}^{O}} \cdot C_{m}^{op} \sum_{k}W_{mk}^{op} \cdot \sum_{t} p_{mkt}
  \end{aligned}
\end{equation}    

\textit{Annualized investment cost method:}
\begin{equation}
 \label{eq:annualized_investment_cost_milestone_weights_proposal1}
  \begin{aligned}
    & \min_{x_{m},p_{mkt}} \quad C^{I} + C^{O}  \\
    & C^{I} = \sum_{m \in \mathcal{M}}\frac{1}{(1+R)^{m}} \cdot x_{m} \cdot \sum_{y = m}^{\min(m+LT-1,Y-1)}\frac{1}{(1+WACC_{m})^{y - m}} \cdot C_{m}^{A} \\ 
    & C^{O} = \sum_{y=0}^{Y-1}\frac{1}{(1+R)^{y}}\sum_{m \in \mathcal{M}} \textcolor{blue}{W_{my}^{O}}\cdot C_{m}^{op} \sum_{k}W_{mk}^{op} \cdot \sum_{t} p_{mkt}
  \end{aligned}
\end{equation} 

Let us first analyze the $C^{I}$ in both equations assuming that the $LT=5$. Therefore, the investment cost in equations in the previous equations becomes:

\textit{Equation (\ref{eq:total_investment_cost_milestone_weights_proposal1}) - total investment cost method:}
\begin{equation*}
    C^{I} = \textcolor{olive}{C_{0}^{T}}x_{0} + \frac{1}{(1+R)^{2}} (C_{2}^{T}-SV_{2})\cdot x_{2} + \frac{1}{(1+R)^{5}} (C_{5}^{T}-SV_{5})\cdot x_{5}
\end{equation*}    

\textit{Equation (\ref{eq:annualized_investment_cost_milestone_weights_proposal1}) - annualized investment cost method:}
\begin{equation*}
  \begin{aligned}
     C^{I} = & \textcolor{olive}{C_{0}^{A}\left(1+\frac{1}{(1+WACC_{0})^{1}}+\frac{1}{(1+WACC_{0})^{2}} +\frac{1}{(1+WACC_{0})^{3}}+\frac{1}{(1+WACC_{0})^{4}}\right) } \cdot x_{0} \\ 
             & +\frac{1}{(1+R)^{2}}C_{2}^{A}\left(1+\frac{1}{(1+WACC_{2})^{1}}+\frac{1}{(1+WACC_{2})^{2}} +\frac{1}{(1+WACC_{2})^{3}} \right) \cdot x_{2} \\
             & + \frac{1}{(1+R)^{5}}C_{5}^{A}\cdot x_{5} \\ 
  \end{aligned}
\end{equation*}    

By simply inspecting the methods used, it is clear that the costs represented in the objective function are the same as before when yearly investment decisions were made. For example, if we compare the first term in both costs (highlighted in \textcolor{olive}{olive}), we can see that the values are equivalent according to equations (\ref{eq:total_cost_from_annualized_cost}) and (\ref{eq:annualized_cost_from_total_cost}). Therefore, this proposal deals with the problem of having an asset whose lifetime ends in the middle of two milestone years for the annualized investment cost method in the $C^{I}$.

Regarding the operational cost $C^{O}$, Table \ref{tab:weight_milestone_year_to_year} shows an example of the values in the parameter $W_{my}^{O}$. Please note that it only impacts the operational cost component in equations (\ref{eq:total_investment_cost_milestone_weights_proposal1}) and (\ref{eq:annualized_investment_cost_milestone_weights_proposal1}). For instance, the operational cost of non-modelled year 1 will be represented as half the values of milestone years 0 and 2. Non-modelled years 3 and 4 are also a linear combination of the adjacent milestone years 2 and 5. This example uses a linear distribution; however, since it is an input parameter, any other type of distribution can be used.

\begin{table}[ht!]
\centering
\begin{tabular}{|c|C{1cm}|C{1cm}|C{1cm}|C{1cm}|C{1cm}|C{1cm}|}
\hline
 & \multicolumn{6}{c|}{\textbf{Years $y$}} \\
\cline{2-7}
\textbf{Milestone Year $m$} & \textbf{0} & \textbf{1} & \textbf{2} & \textbf{3} & \textbf{4} & \textbf{5} \\
\hline
\textbf{0} & 1 & 1/2 & 0 & 0 & 0 & 0 \\
\hline
\textbf{2} & 0 & 1/2 & 1 & 2/3 & 1/3 & 0 \\
\hline
\textbf{5} & 0 & 0 & 0 & 1/3 & 2/3 & 1 \\
\hline
\end{tabular}
\caption{Weight of milestone year $m$ on year $y$ ($W_{my}^{O}$)}
\label{tab:weight_milestone_year_to_year}
\end{table}

The parameter $W_{my}^{O}$ allows mapping of each not modelled year to a milestone year, considering the value of money over time. Using equations (\ref{eq:total_investment_cost_milestone_weights_proposal1}) and (\ref{eq:annualized_investment_cost_milestone_weights_proposal1}) and the values of the example, we get:

\begin{equation*}
  \begin{aligned}
    C^{O} = & C_{0}^{op} \sum_{k}W_{0k}^{op} \cdot \sum_{t} p_{0kt} \\
    & + \frac{1}{(1+R)^{1}} \cdot \left( \frac{1}{2} \cdot \sum_{k}W_{0k}^{op} \cdot \sum_{t} p_{0kt} + \frac{1}{2} \cdot \sum_{k}W_{2k}^{op} \cdot \sum_{t} p_{2kt}\right) \\
    & + \frac{1}{(1+R)^{2}} \cdot \left(\sum_{k}W_{2k}^{op} \cdot \sum_{t} p_{2kt}\right) \\
    & + \frac{1}{(1+R)^{3}} \cdot \left( \frac{2}{3} \cdot \sum_{k}W_{2k}^{op} \cdot \sum_{t} p_{2kt} + \frac{1}{3} \cdot \sum_{k}W_{5k}^{op} \cdot \sum_{t} p_{5kt}\right) \\    
    & + \frac{1}{(1+R)^{4}} \cdot \left( \frac{1}{3} \cdot \sum_{k}W_{2k}^{op} \cdot \sum_{t} p_{2kt} + \frac{2}{3} \cdot \sum_{k}W_{5k}^{op} \cdot \sum_{t} p_{5kt}\right) \\
    & + \frac{1}{(1+R)^{5}} \cdot \left( \sum_{k}W_{5k}^{op} \cdot \sum_{t} p_{5kt}\right)     
  \end{aligned}
\end{equation*} 

Here we can see the main drawback of this proposal, since equations (\ref{eq:total_investment_cost_milestone_weights_proposal1}) and (\ref{eq:annualized_investment_cost_milestone_weights_proposal1}) assume that the operational cost is a linear combination of the milestone years (i.e., $m =$ 0, 2, and 5) for non-modelled years (i.e., $y =$ 1, 3, and 4), then the non-modelled years like years 3 and 4 will consider part of the production costs from the milestone year 5 (which is a future year). There is an improvement compared to the standard approach of having weights in the milestone years shown in Section (\ref{sec:milestone_years_weight}); however, it still might have an overestimation of the operational costs in the objective function. The following proposal aims to overcome this issue.

\subsection{Proposal 2 - Mapping investment milestone years to all years and operational milestone years ($W_{my\mu}^{Y}$ and $W_{m\mu}^{M}$)} \label{sec:proposal_double_weight_operation}
Equations (\ref{eq:total_investment_cost_milestone_weights_proposal1}) and (\ref{eq:annualized_investment_cost_milestone_weights_proposal1}) in proposal 1 solve the issue with the investment cost ($C^{I}$) since the representation of the investment cost is the same in both methods. In addition, it improved the value of money over time representation for the operational cost. Nevertheless, there is still a concern about accounting operational costs from future investment decisions. One solution to that situation is shown in equation (\ref{eq:operational_cost_milestone_weights_proposal2}).

\begin{equation}
 \label{eq:operational_cost_milestone_weights_proposal2}
    C^{O} = \sum_{m \in \mathcal{M}}\sum_{y=0}^{Y-1}\frac{1}{(1+R)^{y}}\sum_{\mu \in \mathcal{M}} \textcolor{blue}{W_{my\mu}^{Y}}\cdot C_{m}^{op} \sum_{k}W_{mk}^{op} \cdot \sum_{t} \textcolor{purple}{W_{m\mu}^{M}} \cdot p_{mkt}
\end{equation} 

Where,
\begin{itemize}
    \item The parameter $W_{m\mu}^{M}$ represents the proportion of production assigned to the investment made in milestone year $m$ on the operational milestone year $\mu$.
    \item The parameter $W_{my\mu}^{Y}$ represents the weight of investment milestone year $m$ on year $y$ considering operational milestone year $\mu$. In other words, how much of the operational production of year $\mu$ does one want to consider for the year $y$ of the investment made in year $m$.
\end{itemize}

To understand the parameter $W_{m\mu}^{M}$, let us first consider in our example that equation (\ref{eq:production_limit}) becomes equation (\ref{eq:production_limit_example}) when taking into account the milestone years $m =$ 0, 2, and 5.

\begin{equation}
 \label{eq:production_limit_example}
 \begin{aligned}
    &p_{0kt} \leq x_{0}   &&               &&& \quad \forall k \forall t \\
    &p_{2kt} \leq x_{0} + && x_{2}         &&& \quad \forall k \forall t \\
    &p_{5kt} \leq         && x_{2} + x_{5} &&& \quad \forall k \forall t \\
 \end{aligned}
\end{equation}

Equation (\ref{eq:production_limit_example}) shows that operational production in year 0 is due entirely to the investment made in year 0, while the operational production in year 2 depends on the investments made in both year 0 and year 2, and the operational production in year 5 depends on the investments made in year 2 and year 5, but not year 0, since the lifetime ($LT$=5).
Multiple assumptions can be made here to determine the parameter $W_{m\mu}^{M}$. A naive approach is to distribute equally the number of investments that limits each operational year. For instance, the operational production of year 2 can be split into halves between the investments made in years 0 and 2. Table (\ref{tab:weight_milestone_to_milestone}) shows the $W_{m\mu}^{M}$ values considering this approach for the example in equation (\ref{eq:production_limit_example}).

\begin{table}[ht!]
\centering
\begin{tabular}{|c|C{1cm}|C{1cm}|C{1cm}|}
\hline
 & \multicolumn{3}{c|}{\textbf{Milestone Years $\mu$}} \\
\cline{2-4}
\textbf{Milestone Year $m$} & \textbf{0} & \textbf{2} & \textbf{5} \\
\hline
\textbf{0} & 1 & 1/2 & 0   \\
\hline
\textbf{2} & 0 & 1/2 & 1/2 \\
\hline
\textbf{5} & 0 & 0  & 1/2  \\
\hline
\end{tabular}
\caption{Weight of \textit{investment} milestone year $m$ on \textit{operational} milestone year $\mu$ (\textcolor{purple}{$W_{m\mu}^{M}$})}
\label{tab:weight_milestone_to_milestone}
\end{table}

Another approach to determining the $W_{m\mu}^{M}$ values could be to consider an efficiency on the capacity factor (or degradation) that changes over the milestone years, such as the investments on milestone years with more efficiency get a higher proportion of the operational production. So, for instance, if years 2 and 5 are twice efficient as the technology invested in year 0, then $W_{0,2}^{M} = 1/3$ and $W_{2,2}^{M} = 2/3$. \textcolor{red}{The assumption here is that all investments will be made. In other words, the distribution among the investment decisions is predefined before the optimization decision. However, if $x_{0}=0$ and $x_{2} \neq 0$, then we are underestimating the cost recovery for year 2. The proposed formulation in Section \ref{sec:proposal_double_index_operation} shows a workaround to avoid this assumption.}

For the parameter $W_{my\mu}^{Y}$, a simple approach distributes it linearly between milestone years $m$ and $\mu$. Table (\ref{tab:weight_milestone_m_y_mu}) shows the values for the example we have been considering before. For instance, the non-modelled year 1 is between milestone years 0 and 2; therefore, $W_{0,1,0}^{Y}=1/2$ and $W_{0,1,2}^{Y}=1/2$, meaning that we represent the operational cost of year 1 with half the value in year 0 and half in year 2. Notice that the values for $W_{0,2,2}^{Y}=W_{0,3,2}^{Y}=W_{0,4,2}^{Y}=1$ because the lifetime of the investment on year 0 ends on year 4; so the operational production of milestone year 5 does not consider the investment on year 0; see equation (\ref{eq:production_limit_example}).

\begin{table}[ht!]
\centering
\begin{tabular}{|c|c|C{1cm}|C{1cm}|C{1cm}|}
\hline
 & & \multicolumn{3}{c|}{\textbf{Milestone Year $\mu$}} \\
\cline{3-5}
\textbf{Milestone Year $m$} & \textbf{Year $y$} & \textbf{0} & \textbf{2} & \textbf{5} \\
\hline
\hline
\textbf{0} & \textbf{0} & 1 & 0 & 0 \\
 \cline{2-5}
 & \textbf{1} & 1/2 & 1/2 & 0 \\
 \cline{2-5}
 & \textbf{2} & 0 & 1 & 0 \\
 \cline{2-5}
 & \textbf{3} & 0 & 1 & 0 \\
 \cline{2-5}
 & \textbf{4} & 0 & 1 & 0 \\
 \cline{2-5}
 & \textbf{5} & 0 & 0 & 0 \\
\hline
\hline
\textbf{2} & \textbf{2} & 0 & 1 & 0 \\
\cline{2-5}
 & \textbf{3} & 0 & 2/3 & 1/3 \\
 \cline{2-5}
 & \textbf{4} & 0 & 1/3 & 2/3 \\
 \cline{2-5}
 & \textbf{5} & 0 & 0 & 1 \\
\hline
\hline
\textbf{5} & \textbf{5} & 0 & 0 & 1 \\
\hline
\end{tabular}
\caption{Weight of milestone year $m$ on year $y$ and milestone year $\mu$ (\textcolor{blue}{$W_{my\mu}^{Y}$})}
\label{tab:weight_milestone_m_y_mu}
\end{table}

Using the information of Table \ref{tab:weight_milestone_to_milestone}, Table \ref{tab:weight_milestone_m_y_mu}, $LT=5$, and the equation (\ref{eq:operational_cost_milestone_weights_proposal2}), the total operational cost will be:

\begin{equation*}
  \begin{aligned}
    C^{O} = &  \sum_{m \in \mathcal{M}}\sum_{y=m}^{\min(m+LT-1,Y-1)}\frac{1}{(1+R)^{y}}\sum_{\mu \in \mathcal{M}} \textcolor{blue}{W_{my\mu}^{Y}}\cdot C_{\mu}^{op} \sum_{k}W_{\mu k}^{op} \cdot \sum_{t} \textcolor{purple}{W_{m\mu}^{M}} \cdot p_{\mu kt} \\
    &  = C_{0}^{op} \sum_{k}W_{0,k}^{op} \cdot \sum_{t} p_{0,kt} \\
    &  + \frac{1}{(1+R)^{1}} \frac{1}{2} C_{0}^{op} \sum_{k}W_{0,k}^{op} \cdot \sum_{t} p_{0,kt} + \frac{1}{(1+R)^{1}} \frac{1}{2} C_{2}^{op} \sum_{k}W_{2,k}^{op} \cdot \sum_{t} \frac{1}{2}p_{2,kt}\\
    &  + \frac{1}{(1+R)^{2}} C_{2}^{op} \sum_{k}W_{2,k}^{op} \cdot \sum_{t}  \frac{1}{2} p_{2,kt} \\ 
    &  + \frac{1}{(1+R)^{3}} C_{2}^{op} \sum_{k}W_{2,k}^{op} \cdot \sum_{t}  \frac{1}{2} p_{2,kt} \\ 
    &  + \frac{1}{(1+R)^{4}} C_{2}^{op} \sum_{k}W_{2,k}^{op} \cdot \sum_{t}  \frac{1}{2} p_{2,kt} \\ 
    &  + \frac{1}{(1+R)^{2}} C_{2}^{op} \sum_{k}W_{2,k}^{op} \cdot \sum_{t}  \frac{1}{2} p_{2,kt} \\ 
    &  + \frac{1}{(1+R)^{3}} \frac{2}{3} C_{2}^{op} \sum_{k}W_{2,k}^{op} \cdot \sum_{t} \frac{1}{2} p_{2,kt} + \frac{1}{(1+R)^{3}} \frac{1}{3} C_{5}^{op} \sum_{k}W_{5,k}^{op} \cdot \sum_{t} \frac{1}{2} p_{5,kt}\\
    &  + \frac{1}{(1+R)^{4}} \frac{1}{3} C_{2}^{op} \sum_{k}W_{2,k}^{op} \cdot \sum_{t} \frac{1}{2} p_{2,kt} + \frac{1}{(1+R)^{4}} \frac{2}{3} C_{5}^{op} \sum_{k}W_{5,k}^{op} \cdot \sum_{t} \frac{1}{2} p_{5,kt}\\    
    &  + \frac{1}{(1+R)^{5}} C_{5}^{op} \sum_{k}W_{5,k}^{op} \cdot \sum_{t} \frac{1}{2} p_{5,kt} \\     
    &  + \frac{1}{(1+R)^{5}} C_{5}^{op} \sum_{k}W_{5,k}^{op} \cdot \sum_{t} \frac{1}{2} p_{5,kt} \\ 
  \end{aligned}
\end{equation*} 

Here the first five terms correspond to the operational cost of investment made in year 0, the following four terms to the investment in milestone year 2, and the last one to milestone year 5. Therefore, it improves the previous proposal in Section \ref{sec:proposal_mapping_m_to_y} by correctly accounting for each investment's operational cost in the milestone year, avoiding overestimating the operational cost. However, as mentioned before, the central assumption in this proposal is the predefined distribution of the operational costs before the optimization, which could potentially lead to an underestimation of the operational cost.

\subsection{Proposal 3 - New operational variable} \label{sec:proposal_double_index_operation}
Although proposals in Sections \ref{sec:proposal_mapping_m_to_y} and \ref{sec:proposal_double_weight_operation} improve the standard formulation in Section \ref{sec:milestone_years_weight}, they either overestimate or underestimate the operational costs in the objective function. The main issue is that consolidating the production data into a single variable per year poses a challenge in determining the production outcome for each investment decision in advance. This is because investment decisions are endogenous, resulting from the optimization problem. Therefore, a more detailed formulation is required to accurately determine the production of each investment decision; see equation (\ref{eq:optimization_problem_new_operational_var}).

\begin{equation}
 \label{eq:optimization_problem_new_operational_var}
 \begin{aligned}
    & \min_{x_{m},p_{m\mu kt},p_{\mu kt}} \quad C^{I} + C^{O} \\
    &s.t. \\
    &p_{m\mu kt} \leq x_{m} \quad \forall m \forall \mu \forall k \forall t | m \le \mu \leq min(m+LT-1,Y-1) \\
    &p_{\mu kt} = \sum_{m \ge \max(\mu-LT+1,0)}^{\mu}p_{m\mu kt} \quad \forall \mu \forall k \forall t
 \end{aligned}
\end{equation}

The new variable $p_{m\mu kt}$ represents the production imputed to the investment in year $m$ for period $t$ and representative $k$ in operational year $\mu$. In the example, it will be:

\begin{equation*}
 \begin{aligned}
   &p_{0,0,k,t} && \leq x_{0} \quad \forall k \forall t \\
   &p_{0,2,k,t} && \leq x_{0} \quad \forall k \forall t \\
   &p_{2,2,k,t} && \leq x_{2} \quad \forall k \forall t \\
   &p_{2,5,k,t} && \leq x_{2} \quad \forall k \forall t \\
   &p_{5,5,k,t} && \leq x_{5} \quad \forall k \forall t \\
   &p_{0,k,t}   &&    = p_{0,0,k,t} \quad  \forall k \forall t \\
   &p_{2,k,t}   &&    = p_{0,2,k,t} + p_{2,2,k,t} \quad  \forall k \forall t \\
   &p_{5,k,t}   &&    = p_{2,5,k,t} + p_{5,5,k,t} \quad  \forall k \forall t \\   
 \end{aligned}
\end{equation*}

Then the operational cost can be easily determined using the new variables and the parameter $W_{my\mu}^{Y}$ as in equation (\ref{eq:operational_cost_milestone_weights_proposal3}).

\begin{equation}
 \label{eq:operational_cost_milestone_weights_proposal3}
    C^{O} = \sum_{m \in \mathcal{M}}\sum_{y=0}^{Y-1}\frac{1}{(1+R)^{y}}\sum_{\mu = m}^{min(m+LT-1,Y-1)} \textcolor{blue}{W_{my\mu}^{Y}}\cdot C_{m}^{op} \sum_{k}W_{mk}^{op} \cdot \sum_{t} p_{m\mu kt}
\end{equation} 

For the example we have been using before, equation (\ref{eq:operational_cost_milestone_weights_proposal3}) becomes:

\begin{equation*}
  \begin{aligned}
    C^{O} &  = C_{0}^{op} \sum_{k}W_{0,k}^{op} \cdot \sum_{t} p_{0,0,k,t} \\
          &  + \frac{1}{(1+R)^{1}} \frac{1}{2} C_{0}^{op} \sum_{k}W_{0,k}^{op} \cdot \sum_{t} p_{0,0,k,t} + \frac{1}{(1+R)^{1}} \frac{1}{2} C_{2}^{op} \sum_{k}W_{2,k}^{op} \cdot \sum_{t} p_{0,2,k,t}\\
          &  + \frac{1}{(1+R)^{2}} C_{2}^{op} \sum_{k}W_{2,k}^{op} \cdot \sum_{t} p_{0,2,k,t} \\ 
          &  + \frac{1}{(1+R)^{3}} C_{2}^{op} \sum_{k}W_{2,k}^{op} \cdot \sum_{t} p_{0,2,k,t} \\ 
          &  + \frac{1}{(1+R)^{4}} C_{2}^{op} \sum_{k}W_{2,k}^{op} \cdot \sum_{t} p_{0,2,k,t} \\ 
          &  + \frac{1}{(1+R)^{2}} C_{2}^{op} \sum_{k}W_{2,k}^{op} \cdot \sum_{t} p_{2,2,k,t} \\ 
          &  + \frac{1}{(1+R)^{3}} \frac{2}{3} C_{2}^{op} \sum_{k}W_{2,k}^{op} \cdot \sum_{t} p_{2,2,k,t} + \frac{1}{(1+R)^{3}} \frac{1}{3} C_{5}^{op} \sum_{k}W_{5,k}^{op} \cdot \sum_{t} p_{2,5,k,t}\\
          &  + \frac{1}{(1+R)^{4}} \frac{1}{3} C_{2}^{op} \sum_{k}W_{2,k}^{op} \cdot \sum_{t} p_{2,2,k,t} + \frac{1}{(1+R)^{4}} \frac{2}{3} C_{5}^{op} \sum_{k}W_{5,k}^{op} \cdot \sum_{t} p_{2,5,k,t}\\    
          &  + \frac{1}{(1+R)^{5}} C_{5}^{op} \sum_{k}W_{5,k}^{op} \cdot \sum_{t} p_{2,5,k,t} \\     
          &  + \frac{1}{(1+R)^{5}} C_{5}^{op} \sum_{k}W_{5,k}^{op} \cdot \sum_{t} p_{5,5,k,t} \\ 
  \end{aligned}
\end{equation*} 

Here it is possible to see that all the issues stated in Section \ref{sec:proposals_main_sec} are solved. In addition, this proposal avoids both overestimating and underestimating the operational costs. However, it does all this at the cost of having more variables in the optimization problem. For instance, there are five more variables in the numerical example. 

\section{Formulations Summary}
Throughout this document, we analyzed various formulations and modifications to accurately account for the investments and operational costs of multi-year investment optimization problems in energy systems. We have demonstrated that both the total investment cost method and the annualized investment cost method lead to the same representation if executed carefully. Additionally, the operational cost representation presents challenges in avoiding issues with the value of money over time and the asset's lifetime. The proposals in Section \ref{sec:proposals_main_sec} provide a comprehensive representation that improves the current standard methods in the literature. Here is a summary of all the discussed formulations:

\begin{itemize}
    \item Weights on milestone years (standard approach)
    \begin{itemize}
        \item Equations: (\ref{eq:total_investment_cost_milestone_weights}) and (\ref{eq:annualized_investment_cost_milestone_weights})
        \item \textcolor{Green}{PROS:} Simple and straightforward to use.
        \item \textcolor{Red}{CONS:} It does not account appropriately for the investment costs in the annualized formulation nor the value of money over time
    \end{itemize}
    \item Proposal 1 - Mapping investment milestone years to all years ($W_{my}^{O}$)
    \begin{itemize}
        \item Equations: (\ref{eq:total_investment_cost_milestone_weights_proposal1}) and (\ref{eq:annualized_investment_cost_milestone_weights_proposal1})
        \item \textcolor{Green}{PROS:} It accounts properly for the investment costs in the annualized formulation as well as the value of money over time
        \item \textcolor{Red}{CONS:} It overestimates the operational costs
    \end{itemize}
    \item Proposal 2 - Mapping investment milestone years to all years and operational milestone years ($W_{my\mu}^{Y}$ and $W_{m\mu}^{M}$)
    \begin{itemize}
        \item Equations: (\ref{eq:total_investment_cost_milestone_weights_proposal1}) and (\ref{eq:annualized_investment_cost_milestone_weights_proposal1}), but changing the $C^{O}$ for the one in (\ref{eq:operational_cost_milestone_weights_proposal2})
        \item \textcolor{Green}{PROS:} It accounts properly for the investment costs in the annualized formulation as well as the value of money over time 
        \item \textcolor{Red}{CONS:} It underestimates the operational costs 
    \end{itemize}
    \item Proposal 3 - New operational variable
    \begin{itemize}
        \item Equations:  (\ref{eq:total_investment_cost_milestone_weights_proposal1}) and (\ref{eq:annualized_investment_cost_milestone_weights_proposal1}), but changing the $C^{O}$ for the one in  (\ref{eq:operational_cost_milestone_weights_proposal3})
        \item \textcolor{Green}{PROS:} It accounts for all the costs in the objective function properly
        \item \textcolor{Red}{CONS:} More variables in the optimization problem
    \end{itemize}    
\end{itemize}

\subsection{Final Comments}
In the event that the investment is a fixed value and the optimization problem is centred solely on operational costs, Proposal 2 can be executed without underestimations on production costs.

By extending the values of weights beyond the last year of modelling in any of the proposals, we can consider new production costs instead of solely relying on the recovery assumption. This approach broadens the scope of economic assumptions in analysing assets barely produced in the last modelled year. Furthermore, it enables us to have a lengthier cash flow projection, which is instrumental in double-checking the profitability of the investment without necessarily increasing the number of milestone years.

Finally, this document's future versions will analyse the formulation changes when considering decommissioning variables and vintage investment representation.

\printbibliography

\end{document}